\documentclass[a4paper]{amsart}

\usepackage{amssymb,amsmath}
\usepackage{latexsym}
\usepackage{eucal}
\makeatletter
\@addtoreset{equation}{section}\makeatother

\newtheorem{theo}{Theorem}[section]
\newtheorem{lem}[theo]{Lemma}

\newtheorem{cor}[theo]{Corollary}

\newcommand{\mc}{\mathcal}
\newcommand{\rr}{\mathbb{R}}
\newcommand{\nn}{\mathbb{N}}
\newcommand{\cc}{\mathbb{C}}
\newcommand{\hh}{\mathbb{H}}
\newcommand{\zz}{\mathbb{Z}}

\newcommand{\la}{\lambda}
\newcommand{\eps}{\epsilon}

\newcommand{\pl}{\partial}
\newcommand{\x}{\times}

\newcommand{\til}{\widetilde}
\newcommand{\bbar}{\overline}

\newcommand{\ddens}{\Gamma_{0}^{\demi}}

\newcommand{\cjd}{\rangle}
\newcommand{\cjg}{\langle}

\newcommand{\demi}{\frac{1}{2}}
\newcommand{\ndemi}{\frac{n}{2}}
\newcommand{\tra}{\textrm{Tr}}
\newcommand{\beq}{\begin{equation}}
\newcommand{\eeq}{\end{equation}}

\newcommand{\rang}{\textrm{rank}}
\newcommand{\rangp}{\textrm{Rank}}
\newcommand{\trans}{{^t}\!}
\newcommand{\indic}{\operatorname{1\negthinspace l}}
\def\qed{\hfill$\square$}
\begin{document}

\title{Resonances and scattering poles
on asymptotically hyperbolic manifolds}
\author[Colin Guillarmou]{Colin Guillarmou}
\address{Laboratoire de Math\'ematiques Jean Leray\\
         UMR 6629 CNRS/Universit\'e de Nantes \\
         2, rue de la Houssini\`ere \\
         BP 92208 \\
         44322 Nantes Cedex 03\\
         France}
     \email{cguillar@math.univ-nantes.fr}
\subjclass[2000]{Primary 58J50, Secondary 35P25}
\maketitle

\begin{abstract}
\noindent On an asymptotically hyperbolic manifold $(X,g)$,
we show that the poles (called resonances)
of the meromorphic extension of the resolvent
$(\Delta_g-\la(n-\la))^{-1}$
coincide, with multiplicities, with the poles (called scattering poles)
of the renormalized scattering operator, except
for the points of $\ndemi-\nn$.
At each $\la_k:=\ndemi-k$ with $k\in\nn$,
the resonance multiplicity $m(\la_k)$ and the scattering pole
multiplicity $\nu(\la_k)$ do not always coincide:
$\nu(\la_k)-m(\la_k)$ is the dimension of the kernel
of a differential operator on the boundary $\pl\bar{X}$
introduced by Graham and Zworski;
in the asymptotically Einstein case, this operator is
the k-th conformal Laplacian.
\end{abstract}
\vspace{0.5cm}

\section{Introduction}

The purpose of this work is to give a `more direct' proof of the result
of Borthwick and Perry \cite{PB} about the equivalence between
resolvent resonances and scattering poles, notably in order
to analyze the special points $(\frac{n-k}{2})_{k\in\nn}$
that they did not deal with. This problem is especially interresting
on convex co-compact hyperbolic quotients since these are the
scattering poles (not the resonances) which appear in the divisor
of Selberg's zeta function associated to the group
(cf. Patterson-Perry \cite{PP}).

Let $\bar{X}=X\cup\pl\bar{X}$ a $n+1$-dimensional smooth compact
manifold with boundary and $x$ a defining function for the
boundary, that is a smooth function $x$ on $\bar{X}$ such that
\[x\geq 0,\quad \pl\bar{X}=\{m\in\bar{X}, x(m)=0\},\quad
dx|_{\pl\bar{X}}\not=0\]
We say that a smooth metric $g$ on the
interior $X$ of $\bar{X}$ is \textsl{conformally compact} if
$x^2g$ extends smoothly as a metric to $\bar{X}$. An
\textsl{asymptotically hyperbolic manifold} is a conformally
compact manifold such that for all $y\in\pl\bar{X}$, all sectional
curvatures at $m\in X$ converge to $-1$ as $m\to y$.
Notice that convex co-compact hyperbolic quotients are included in this
class of manifolds.
An asymptotically hyperbolic manifold is necessarily
complete and the spectrum of its Laplacian
$\Delta_g$ acting on functions consists of absolutely continuous
spectrum $[\frac{n^2}{4},\infty)$ and a finite set of eigenvalues
$\sigma_{pp}(\Delta_g)\subset (0,\frac{n^2}{4})$. The resolvent
$(\Delta_g-z)^{-1}$ is a meromorphic family on $\cc\setminus
[\frac{n^2}{4},\infty)$ of bounded operators and the new parameter
$z=\la(n-\la)$ with $\Re(\la)>\ndemi$ induces a modified resolvent
\[R(\la):=(\Delta_g-\la(n-\la))^{-1}\]
which is meromorphic on $\{\Re(\la)>\ndemi\}$, its poles being the
points $\la_e$ such that $\la_e(n-\la_e)\in\sigma_{pp}(\Delta_g)$.
Mazzeo and Melrose \cite{MM} have constructed the
finite-meromorphic extension (i.e. with poles whose residue
is a finite rank operator) of $R(\la)$ on $\cc\setminus \demi(n-\nn)$. We
proved in a previous work \cite{GU} that this extension is
finite-meromorphic on $\cc$ if and only if  the metric is even
in the sense that there exists a boundary defining function $x$
such that the metric can be expressed by
\beq\label{even}
g=\frac{dx^2+h(x^2,y,dy)}{x^2}
\eeq
in the collar $[0,\eps)\x\pl\bar{X}$ induced by $x$,
with $h(z,y,dy)$ smooth up to $\{z=0\}$.
We will only consider these cases
of even metrics to simplify the statements, but our result works as
long as the studied singularity is a pole of
finite multiplicity for the resolvent.\\

The poles of the extension $R(\la)$ are called \textsl{resonances}
and the multiplicity of a resonance $\la_0$ is defined by
\[m(\la_0):=\rang
\int_{C(\la_0,\eps)}(n-2\la)R(\la)d\la=\rang\textrm{Res}_{\la_0}((n-2\la)R(\la))\]
where $C(\la_0,\eps)$ is a circle around $\la_0$ with radius
$\eps>0$ chosen sufficiently small to avoid other resonances in
$D(\la_0,\eps)$ and $\textrm{Res}$ means the residue.
In other words, this is the rank of the residue
at $z_0=\la_0(n-\la_0)$ of the resolvent as a function of
$z=\la(n-\la)$.\\

The scattering operator $S(\la)$ is the operator on $\pl\bar{X}$ 
defined as follows: let
$\la\in \{\Re(\la)=\ndemi\}$ and $\la\not=\ndemi$, for all $f_0\in
C^\infty(\pl\bar{X})$ there exists a unique solution $F(\la)$ of
the problem
\[(\Delta_g-\la(n-\la))F(\la)=0,\quad F(\la)=x^\la f_-+x^{n-\la}f_+\]
\[f_-,f_+ \in C^\infty(\bar{X}), \quad f_+|_{\pl\bar{X}}=f_0\]
we then set $S(\la)$ the operator $S(\la): f_0\to
f_-|_{\pl\bar{X}}$. In fact we should use half-densities and
define $S(\la)$ on conormal bundles on $\pl\bar{X}$ to get
invariance with respect to $x$, but this is dropped here.
Joshi and S\'a Barreto showed \cite{JSB} that
this family of operators extends meromorphically in $\cc\setminus
\demi(n-\nn)$ in the sense of pseudo-differential operators on
$\pl\bar{X}$ and that $S(\la)$ has the principal symbol
\beq\label{symboleprincipal} \sigma_{0}\left (S(\la)\right
)=c(\la)\sigma_{0}\left (\Lambda^{2\la-n}\right ),\textrm{ with } \Lambda:=(1+\Delta_{h_0})^\demi, \quad
c(\la):=2^{n-2\la}\frac{\Gamma(\ndemi-\la)}{\Gamma(\la-\ndemi)}\eeq
and $h_0:=x^2g|_{T\pl\bar{X}}$, which leads to the factorization 
(see \cite{P1,GZ3,PP,PB} for a similar approach)
\beq\label{stilde}
\til{S}(\la):=c(n-\la)\Lambda^{-\la+\ndemi}S(\la)\Lambda^{-\la+\ndemi}=
1+K(\la) \eeq
with $K(\la)$ compact finite-meromorphic. It is clear that
the poles of $S(\la)$ and $\til{S}(\la)$ coincide except
for the points of $\ndemi+\zz$. A pole $\la_0$ of
$\til{S}(\la)$ is called a \textit{scattering pole} and we define
its multiplicity by
\[\nu(\la_0):=-\tra\left(\frac{1}{2\pi i}\int_{C(\la_0,\eps)}
\til{S}'(\la)\til{S}^{-1}(\la)d\la\right)=-\tra\textrm{Res}_{\la_0}
(\til{S}'(\la)\til{S}^{-1}(\la)).\]
Using a method close to that of Guillop\'e-Zworski \cite{GZ3}
and Gohberg-Sigal theory \cite{GS}, we then obtain the

\begin{theo}\label{multiplicite}
Let $(X,g)$ be an asymptotically hyperbolic manifold with $g$ even 
in the sense of (\ref{even}) and let $\la_0\in\{\Re(\la)<\ndemi\}$ such that
$\la_0\notin \left\{\la\in\cc; \la(n-\la)\in\sigma_{pp}(\Delta_g)\right\}
\cap \demi(n-\nn)$.
Then $\la_0$ is a pole of $R(\la)$ if and only if it is a
pole of $S(\la)$ and we have
\beq\label{egalitemult}
m(\la_0)=m(n-\la_0)+\nu(\la_0)-\indic_{\ndemi-\nn}(\la_0)\dim\ker
\textrm{Res}_{n-\la_0}S(\la)
\eeq
where $\indic_{\ndemi-\nn}$ is the characteristic function
of $\ndemi-\nn$ and $\textrm{Res}$ means the residue.
\end{theo}

\textsl{Remark 1}: the term $m(n-\la_0)$ vanishes when $\la_0(n-\la_0)\notin
\sigma_{pp}(\Delta_g)$ and that (\ref{egalitemult})
can be extended to the line
$\{\Re(\la)=\ndemi\}$ by using that $R(\la)$ and $\til{S}(\la)$
are continuous on this line except possibly at $\ndemi$,
where only $R(\la)$ can have a pole; in this case
$\nu(\la_0)=0$ and (\ref{egalitemult}) is satisfied.\\

\textsl{Remark 2}: the additional term introduced  at $\la_0=\ndemi-k$
is exactly the dimension of the kernel of the operator $p_{2k}$ defined by
Graham-Zworski in \cite[Prop. 3.5]{GRZ}. Therefore it
only depends on the
$2k$ first derivatives of the metric at the boundary.
When the manifold is asymptotically Einstein, this is
\[\dim\ker \textrm{Res}_{\ndemi+k}S(\la)=\dim\ker P_k\]
$P_k$ being the k-th conformally invariant power of the Laplacian
(cf. \cite{GRZ}), which depends only on the conformal class of the metric
$h_0=x^2g|_{T\pl\bar{X}}$ at the boundary. If $n$ is even, it is worth 
noting that $\dim\ker p_n\geq 1$ since $p_n$ always annihilates 
constants. Moreover, if $(\pl\bar{X},h_0)$ is conformally flat with $(X,g)$
asymptotically Einstein, the additional term is
$\dim\ker P_k=H_0(\pl\bar{X})$, the number of connected
components of the boundary.

The recent formula obtained by Patterson-Perry \cite{PP} and Bunke-Olbrich \cite{BO}
for the divisor at $\la_0\in\cc$
of Selberg's zeta function on a convex co-compact hyperbolic
quotient always makes the `spectral term' $\nu(\la_0)$ appear
and an additional `topological term' (an integer multiple of
the Euler characteristic) comes when $\la_0\in-\nn_0$.
As a matter of fact, the `spectral term' at $\la_0=\ndemi-k$ (with $k\in\nn$)
could be splitted in a `resonance term' $m(\la_0)$ and a
`conformal term' $\dim\ker p_{2k}$ with $p_{2k}$ the residue
of $S(\la)$ at $\ndemi+k$. Notice also that for $\la_0\in\ndemi-\nn$,
$m(\la_0)$ can be $0$ though $\nu(\la_0)$ is not
(this is the case of $\hh^{n+1}$ when $n+1$ is odd).

Moreover the Poisson formula obtained by
Perry \cite{P2} for convex co-compact quotients
is used to give a lower bound of poles
of $\til{S}(\la)$ (with multiplicity $\nu(\la_0)$)
in a disc $D(\ndemi,R)\subset\cc$ with radius $R$.
It is clear that the number of these poles is bigger than
the number of resonances, in view of
Theorem \ref{multiplicite}. In the trivial case of
$\hh^{n+1}$ with $n+1$ odd, we notably
have no resonance though the number of poles of
$\til{S}(\la)$ in $D(\ndemi,R)$ is $CR^{n+1}$.
However, in dimension $n+1=2$, the explicit formula of the scattering matrix
for a hyperbolic funnel by Guillop\'e-Zworski \cite{GZ1} or
the work of Bunke-Olbrich \cite[Prop.4.3]{BO1}
show that the conformal term cancels, so
$\nu(\la_0)=m(\la_0)$ (modulo the discrete spectrum).

To conclude it would be interresting to study
the dimension of the kernels of the conformal Laplacians on
such quotients to use Perry's results and give a
lower bound of the number of resonances in a disc.\\

\textbf{Acknowledgements.} I would like to thank L. Guillop\'e for
help and comments. I also thank R. Graham, M. Olbrich,
M. Zworski and G.Vodev for useful discussions.

\section{Background on multiplicities}

Let $\mc{H}_1$, $\mc{H}_2$ some Hilbert spaces.
If $M(\la)$ is meromorphic on an open set $U\subset \cc$
with values in the space $\mc{L}(\mc{H}_1,\mc{H}_2)$ of bounded linear operators
and if $\la_0$ is a pole of $M(\la)$,
there exists a neighborhood $V_{\la_0}$ of $\la_0$, an integer $p>0$ and
some $(M_i)_{i=1,\dots,p}$ in $\mc{L}(\mc{H}_1,\mc{H}_2)$ such that
for $\la\in V_{\la_0}\setminus \{\la_0\}$
\beq\label{serielaurent}
M(\la)=\Xi_{\la_0} (M(\la))+H(\la),
\eeq
\[\Xi_{\la_0}(M(\la))=\sum_{i=1}^pM_i(\la-\la_0)^{-i}, \quad
H(\la)\in\mc{H}ol(V_{\la_0},\mc{L}(\mc{H}_1,\mc{H}_2)).\]
We will call $\Xi_{\la_0}(M(\la))$ the polar part of $M(\la)$ at $\la_0$,
$p$ the order of the pole $\la_0$,
$M_1=\textrm{Res}_{\la_0}M(\la)$ the residue of $M(\la)$ at $\la_0$,
$m_{\la_0}(M(\la)):=\rang M_1$ the multiplicity of $\la_0$ and
\[\rangp_{\la_0}M(\la) := \dim\sum_{i=1}^p \textrm{Im}(M_i)\]
the total polar rank of $M(\la)$ at $\la_0$.
Finally, a meromorphic family of operators in $\mc{L}(\mc{H}_1,\mc{H}_2)$
whose poles have finite total polar rank will be called
finite-meromorphic.\\

Assume now that $\mc{H}_1=\mc{H}_2$;
taking essentially Gohberg-Sigal notations \cite{GS},
a root function of $M(\la)$ at $\la_0$ is a function
$\varphi(\la)\in\mc{H}ol(V_{\la_0},\mc{H}_1)$ such that
$\lim_{\la\to\la_0}M(\la)\varphi(\la)=0$ and
$\varphi(\la_0)\not=0$, the vanishing order
of $M(\la)\varphi(\la)$ being called the multiplicity of
$\varphi(\la)$. The vector $\varphi_0:=\varphi(\la_0)$ is
called an eigenvector of $M(\la)$ at $\la_0$
and the set of eigenvectors of $M(\la)$
at $\la_0$ form a vectorial subspace of $\mc{H}_1$
denoted $\ker_{\la_0} M(\la)$. The rank of an eigenvector $\varphi_0$ is
defined as being the supremum of the multiplicities
of the root functions $\varphi(\la)$ of $M(\la)$ at $\la_0$ such that
$\varphi(\la_0)=\varphi_0$.
If $\dim\ker_{\la_0} M(\la)=\alpha<\infty$ and
the ranks of all eigenvectors are finite,
a canonical system of eigenvectors is a basis
$(\varphi^{(i)}_{0})_{i=1,\dots,\alpha}$ of $\ker_{\la_0} M(\la)$ such
that the ranks of $\varphi^{(i)}_0$ have the following property: the rank
of $\varphi^{(1)}_0$ is the maximum of the ranks of all
 eigenvectors of $M(\la)$
at $\la_0$ and the rank of
$\varphi^{(i)}_0$ is the maximum of the ranks of all eigenvectors
in a direct complement of $\textrm{Vect}(\varphi^{(1)}_0,
\dots,\varphi^{(i-1)}_0)$ in $\ker_{\la_0}M(\la)$.
A canonical system of eigenvectors is not unique but the family of
ranks of its eigenvectors does not depend on the choice of
the canonical system.
We then denote $r_i=\varphi^{(i)}_0$ the partial null multiplicities of
$M(\la)$ at $\la_0$ and
\[N_{\la_0}(M(\la))=\sum_{i=1}^\alpha r_i\]
the null multiplicity of $M(\la)$ at $\la_0$.\\

Assume that $M(\la)$ is meromorphic family of Fredholm operators
in $\mc{L}(\mc{H}_1)$ and $\la_0$ a pole of finite total polar
rank. If the index of $(M(\la)-\Xi_{\la_0}(M(\la))|_{\la=\la_0}$
is $0$, Gohberg and Sigal \cite{GS} show that there exist some
holomorphically invertible operators $U_1(\la)$ and $U_2(\la)$
near $\la_0$, some orthogonal projections $(P_l)_{l=0,\dots,m}$
and some non zero integers $(k_l)_{l=1,\dots,m}$ such that
\beq\label{factorisation}
M(\la)=U_1(\la)\left(P_0+\sum_{l=1}^m(\la-\la_0)^{k_l}P_l\right)U_2(\la),
\eeq
\[P_iP_j=\delta_{ij}P_j, \quad \rang (P_l)=1 \textrm{ for } l=1,\dots,m,
\quad \dim (1-P_0)<\infty.
\]
If moreover $M(\la)$ has a meromorphic inverse $M^{-1}(\la)$ (ie.
when $P_0+\sum_{l=1}^mP_l=1$) then $\la_0$ is at most a pole
of finite total polar rank of $M^{-1}(\la)$ and
\beq\label{factorisation2}
M^{-1}(\la)=U_2^{-1}(\la)\left(P_0+\sum_{l=1}^m(\la-\la_0)^{-k_l}
P_l\right)U_1^{-1}(\la). \eeq
It is important to notice that the set of partial null multiplicities
remains invariant under multiplication by a holomophically invertible
family of operators (cf. \cite{GS}).
In view of (\ref{factorisation}) and (\ref{factorisation2}),
it is then easy to see that
\[\dim\ker_{\la_0} M(\la)=\sharp \{l; k_l>0\}, \quad
\dim\ker_{\la_0} M^{-1}(\la)=\sharp \{l; k_l<0\}\]
and that the set of partial null multiplicities of
$M(\la)$ (resp. $M^{-1}(\la)$) at $\la_0$ is $\{k_l; k_l>0\}$
(resp. $\{k_l;k_l<0\}$). We deduce
\[N_{\la_0}(M(\la))=\sum_{k_l>0}k_l, \quad N_{\la_0}(M^{-1}(\la))=\sum_{k_l<0}k_l\]
and from the factorization (\ref{factorisation}) Gohberg-Sigal \cite{GS} obtain
the generalized logarithmic residue theorem
\beq\label{gohbergsigal}
\tra\left(
\textrm{Res}_{\la_0}(M'(\la)M^{-1}(\la))\right)=N_{\la_0}(M(\la))
-N_{\la_0}(M^{-1}(\la)).
\eeq
This integer is essentially the order of the zero or the pole 
of $\det(M(\la))$ at $\la_0$
(when $det(M(\la))$ exists).\\

To conclude, let $M(\la)$ be a meromorphic family of
Fredholm operators with index $0$ in
$\mc{L}(\mc{H}_1)$ and $\la_0$ a pole of finite total polar rank.
We write $M(\la)$ as in (\ref{factorisation}) and 
if $L(\la):=(\la-\la_0)^{-1}M(\la)$, we deduce that 
$\dim\ker_{\la_0} L(\la)=\sharp\{l;k_l>1\}$, the set
of partial null multiplicities
of $L(\la)$ at $\la_0$ is $\{k_l-1;k_l>1\}$ and
\beq\label{ndela0}
N_{\la_0}(L(\la))=\sum_{k_l>1}(k_l-1)=\sum_{k_l>0}(k_l-1)=
N_{\la_0}(M(\la))-\dim\ker_{\la_0}M(\la).
\eeq
This formula will be essential for what follows
since the scattering operator $S(\la)$ is not
finite-meromorphic near $\ndemi+k$ (with $k\in\nn)$ whereas
$(\la-\ndemi-k)S(\la)$ is.

\section{Resonances and scattering poles}

\subsection{Stretched products, half-densities}

To begin, let us introduce a few notations and recall some basic
things on stretched products and singular half-densities
(the reader can refer to Mazzeo-Melrose \cite{MM}, Melrose \cite{M1}
for details). Let $\bar{X}$ a smooth compact manifold with boundary and $x$
a boundary defining function. The
manifold $\bar{X}\x\bar{X}$ is a smooth manifold with corners,
whose boundary hypersurfaces are diffeomorphic to
$\pl\bar{X}\x\bar{X}$ and $\bar{X}\x\pl\bar{X}$, and defined by
the functions $\pi_L^*x$, $\pi_R^*x$ ($\pi_L$ and $\pi_R$ being
the left and right projections from $\bar{X}\x\bar{X}$ onto
$\bar{X}$). For notational simplicity, we now write $x$,$x'$
instead of $\pi_L^*x$, $\pi_R^*x$ and let
\[\delta_{\pl\bar{X}}:=\{(m,m)\in\pl\bar{X}\x\pl\bar{X}; m\in \pl\bar{X}\}.\]
The blow-up of $\bar{X}\x\bar{X}$ along the diagonal $\delta_{\pl\bar{X}}$ of
$\pl\bar{X}\x\pl\bar{X}$ will be noted $\bar{X}\x_0\bar{X}$ and
the blow-down map
\[\beta:\bar{X}\x_0\bar{X}\to \bar{X}\x\bar{X}\]
This manifold with corners has three boundary hypersurfaces
$\mc{T},\mc{B},\mc{F}$ defined by some functions $\rho,\rho',R$
such that $\beta^*(x)=R\rho$, $\beta^*(x')=R\rho'$. Globally,
$\delta_{\pl\bar{X}}$ is replaced by a larger manifold, namely by
its doubly inward-pointing spherical normal bundle of
$\delta_{\pl\bar{X}}$, whose each fiber is a quarter of sphere.
From local coordinates $(x,y,x',y')$ on $\bar{X}\x\bar{X}$, this
amounts to introducing polar coordinates
$(R,\rho,\rho',\omega,y)$ around $\delta_{\pl\bar{X}}$:
\[R:=(x^2+x'^2+|y-y'|^2)^\demi, \quad (\rho,\rho',\omega):=
\left(\frac{x}{R},\frac{x'}{R},\frac{y-y'}{R}\right)\] with
$R,\rho,\rho'\in [0,\infty)$. In these polar coordinates
the Schwartz kernel of $R(\la)$ has a better description.\\

Using evident identifications induced by the inclusions
\[\delta_{\pl\bar{X}}\subset \pl\bar{X}\x\pl\bar{X}\subset
\pl\bar{X}\x\bar{X}\subset\bar{X}\x\bar{X},\]
we denote by $\pl\bar{X}\x_0\bar{X}$ the blow-up of
$\pl \bar{X}\x\bar{X}$ along $\delta_{\pl\bar{X}}$ and
$\pl\bar{X}\x_0\pl\bar{X}$ the blow-up of
$\pl\bar{X}\x\pl\bar{X}$ along $\delta_{\pl\bar{X}}$.
$\til{\beta}$ and $\beta_{\pl}$ are the associated blow-down map
\[\til{\beta}:\pl\bar{X}\x_0\bar{X}\to\pl\bar{X}\x\bar{X},\quad
\beta_{\pl}: \pl\bar{X}\x_0\pl\bar{X}\to\pl\bar{X}\x\pl\bar{X}\]
with $\til{\beta}=\beta|_{\mc{T}}$ and
$\beta_{\pl}=\beta|_{\mc{B}\cap\mc{T}}$.
Note that $r:=R|_{\mc{B}\cap\mc{T}}$ is a defining function of
the boundary of $\pl\bar{X}\x_0\pl\bar{X}$ (which is the lift of
$\delta_{\pl\bar{X}}$ under $\beta_\pl$).\\

Let $\Gamma_0^\demi(\bar{X})$ the line bundle of singular
half-densities on $\bar{X}$, trivialized by
$\nu:=|dvol_g|^{\demi}$, and $\Gamma^\demi(\pl\bar{X})$ the bundle
of half densities on $\pl\bar{X}$, trivialized by
$\nu_0:=|dvol_{h_0}|^\demi$ (where $h_0=x^2g|_{T\pl\bar{X}}$).
From these bundles, one can construct the bundles
$\Gamma_0^\demi(\bar{X}\x\bar{X})$,
$\Gamma_0^\demi(\pl\bar{X}\x\bar{X})$ and
$\Gamma^\demi(\pl\bar{X}\x\pl\bar{X})$ by tensor products and the bundles
$\Gamma_0^\demi(\bar{X}\x_0\bar{X})$,
$\Gamma_0^\demi(\pl\bar{X}\x_0\bar{X})$ and
$\Gamma^\demi(\pl\bar{X}\x_0\pl\bar{X})$ by 
lifting under $\beta$, $\til{\beta}$ and $\beta_{\pl}$ the three
previous bundles. If $M$ denotes $\bar{X}$, $\bar{X}\x\bar{X}$ or
$\pl\bar{X}\x\bar{X}$, we write
$\dot{C}^\infty(M,\Gamma_0^\demi)$ the space of smooths
sections of $\Gamma_0^\demi(M)$ that vanish to all order at all
the boundary hypersurfaces of $M$, and
$C^{-\infty}(M,\Gamma_0^\demi)$ is its topological dual. 
The Hilbert space $L^2(\bar{X},\Gamma_0^\demi)$ and
$L^2(\pl\bar{X},\Gamma^\demi)$ are isomorphic to $L^2(X,dvol_g)$ and 
$L^2(\pl\bar{X},dvol_{h_0})$, they will be denoted $L^2(X)$, 
$L^2(\pl\bar{X})$.\\

For $\alpha\in\rr$, let
$x^{\alpha}L^2(X):=\{f\in C^{-\infty}(\bar{X},\ddens);x^{-\alpha}f\in L^2(X)\}$
and we set $\cjg .,.\cjd$ the symmetric non-degenerate products
\[\cjg u,v\cjd:=\int_{X}uv \textrm{  on } L^2(X), 
\quad \cjg u,v\cjd:=\int_{\pl\bar{X}}uv \textrm{  on }L^2(\pl\bar{X}).\]
We can check by using the first pairing that the
dual space of $x^\alpha L^2(X)$ is isomorphic to $x^{-\alpha}L^2(X)$.
We shall also use the following tensorial notation for $E=x^{\alpha} L^2(X)$
(resp. $E=L^2(\pl\bar{X})$), $\psi,\phi\in E'$
\[\phi\otimes\psi : \left\{
\begin{array}{rcl}
E & \to & E'\\
f&\to& \phi\cjg \psi,f\cjd
\end{array}\right..\]

\subsection{Resolvent}
From \cite{MM,GU}, we know that on an
asymptotically hyperbolic manifold $(X,g)$ with
$g$ even, the modified resolvent
\[R(\la):=(\Delta_g-\la(n-\la))^{-1}\]
extends for all $N>0$ to a finite-meromorphic family of operators
in $\{\Re(\la)>\ndemi-N\}$ with values in
$\mc{L}(x^NL^2(X),x^{-N}L^2(X))$, whose poles, the resonances,
form a discrete set $\mc{R}$ in $\cc$.
Moreover $R(\la)$ is a continuous
operator from $\dot{C}^\infty(\bar{X},\Gamma_0^\demi)$ to
$C^{-\infty}(\bar{X},\Gamma_0^\demi)$, its associated Schwartz
kernel being
\[r(\la)=r_0(\la)+r_1(\la)+r_2(\la)\in C^{-\infty}(\bar{X}\x \bar{X},
\Gamma_0^\demi)\]
with (see \cite{MM} or \cite[Th. 2.1]{PB}):
\[\beta^*(r_0(\la))\in I^{-2}(\bar{X}\x_0\bar{X},\Gamma_0^\demi),\]
\beq\label{k2}
\beta^*(r_1(\la))\in \rho^\la\rho'^\la C^\infty
(\bar{X}\x_0\bar{X},\Gamma_0^\demi), \quad r_2(\la)\in
x^{\la}x'^{\la}C^\infty(\bar{X}\x\bar{X}, \Gamma_0^\demi),
\eeq
where $I^{-2}(\bar{X}\x_0\bar{X},\Gamma_0^\demi)$ denotes the set of
conormal distributions of order $-2$ on $\bar{X}\x_0\bar{X}$
associated to the closure of the lifted interior diagonal
\[\bbar{\beta^{-1}(\{(m,m)\in\bar{X}\x\bar{X}; m\in X\})}\]
and vanishing to infinite order at $\mc{B}\cup\mc{T}$ (note that
the lifted interior diagonal only intersects the topological
boundary of $\bar{X}\x_0\bar{X}$ at $\mc{F}$,
and it does transversally).
Moreover, $(\rho\rho')^{-\la}\beta^{*}(r_1(\la))$ and $(xx')^{-\la}r_2(\la)$
are meromorphic in $\la\in\cc $ and $r_0(\la)$ is the
kernel of a holomorphic family of operators
\[R_0(\la)\in \mc{H}ol(\cc, \mc{L}
(x^{\alpha}L^2(X),x^{-\alpha}L^2(X))),\quad \forall \alpha\geq 0.\]
Note also that Patterson-Perry arguments \cite[Lem.4.9]{PP}
prove that $R(\la)$ does not have poles on the line
$\{\Re(\la)=\ndemi\}$, except maybe $\la=\ndemi$. The set of poles
of $R(\la)$ in the half plane $\{\Re(\la)>\ndemi\}$ is
$\{\la_e;\Re(\la_e)>\ndemi, \la_e(n-\la_e)\in
\sigma_{pp}(\Delta_g)\}$,
they are first order poles and their residue is  
\beq\label{residuvp}
\textrm{Res}_{\la_e}R(\la)=(2\la_e-n)^{-1}\sum_{k=1}^p \phi_k\otimes
\phi_k, \quad \phi_k\in x^{\la_e}C^\infty(\bar{X},\Gamma_0^\demi),
\eeq
where $(\phi_k)_{k=1,\dots,p}$ are the normalized eigenfunctions of
$\Delta_g$ for the eigenvalue $\la_e(n-\la_e)$.
One can see by a Taylor expansion at $x=0$ of the eigenvector
equation that if $x^{-\la_e+\ndemi}\phi_k|_{\pl\bar{X}}=0$ then
$\phi_k\in \dot{C}^\infty(\bar{X},\Gamma_0^\demi)$, which is excluded
according to Mazzeo's results \cite{MA}.\\

To simplify the notations, we shall set
$z(\la):=\la(n-\la)$ the holomorphically invertible
function from $\Re(\la)<\ndemi$ to $\cc\setminus
[\frac{n^2}{4},\infty)$.\\

For the poles of $R(\la)$ in $\{\Re(\la)<\ndemi\}$,
we use Lemma 2.4 and 2.11 of \cite{GZ3} to show the

\begin{lem}\label{partiesinguliere}
Let $\la_0\in \mc{R}$ and $N$ such that $\ndemi>\Re(\la_0)>\ndemi-N$, then
in a neighbourhood $V_{\la_0}$ of $\la_0$ we have the decomposition
\beq\label{factorisationderla}
R(\la)=\trans\Phi F_1(\la)\left(\sum_{j=1}^m(z(\la)-z(\la_0))^{k_j}
P_j\right)F_2(\la)\Phi+H(\la),
\eeq
with $m\in\nn$, $k_1,\dots k_m\in-\nn$,
\[H(\la)\in\mc{H}ol(V_{\la_0},\mc{L}(x^NL^2(X),x^{-N}L^2(X))),\quad
F_i(\la)\in \mc{H}ol(V_{\la_0},\mc{L}(\cc^q)),\]
where $q=-\sum_{j=1}^m k_j=m_{\la_0}(z'(\la)R(\la))$ is the multiplicity
of the resonance $\la_0$,
$(P_j)_{j=1,\dots,m}$ are some orthogonal projections on $\cc^q$
such that $P_iP_j=\delta_{ij}P_j$ and $\rang (P_j)=1$, $\Phi$ is defined by
\[\Phi :\left\{
\begin{array}{rcl}
x^NL^2(X) & \to & \cc^q\\
f & \to & (\cjg \psi_l,f\cjd)_{l=1,\dots,q}
\end{array}\right.,\]
$(\psi_l)_{l=1,\dots,q}$ being a basis of $\textrm{Im}(A)$ with
$A:=\textrm{Res}_{\la_0}(z'(\la)R(\la))$. Moreover we have
\beq\label{decomp}
\textrm{Im}(A)\subset \sum_{j=0}^{p-1}x^{\la_0}\log^j(x) C^\infty(\bar{X},\Gamma_0^\demi)
\eeq
with $p$ the order of the pole $\la_0$ of $R(\la)$.
\end{lem}
\textsl{Proof}: it suffices to use Lemmas 2.4 and 2.11 of \cite{GZ3}
but we factorize the resolvent and not the scattering operator.
The arguments used in these lemmas are essentially that the polar part
of $R(\la)$ be expressed by
\[\Xi_{\la_0}(R(\la))=\Xi_{\la_0}\left(\sum_{i=1}^p\frac{(\Delta_g-z(\la_0))^{i-1}A}
{(z(\la)-z(\la_0))^i}\right)\]
and the factorization into its Jordan form of the nilpotent matrix
of $\Delta_g-z(\la_0)$ acting on $\textrm{Im}(A)$.
Observe that the elliptic regularity implies that the elements
of $\textrm{Im}(A)$ are smooth in $X$.

To study the structure of the Schwartz kernel $a_j$ of $A_j$, we
first consider the following operator
\beq\label{rtilde}
\til{R}(\la):=x^{-\la+\ndemi}R(\la)x^{-\la+\ndemi}
\eeq
in a disc $D(\la_0,\eps)$ around $\la_0$ with radius $\eps$. If $\eps$ is
taken sufficiently small, $\til{R}(\la)$ is meromorphic in this
disc with values in $\mc{L}(x^{2\eps} L^2(X),x^{-2\eps}L^2(X))$,
$\la_0$ is the only pole and its order is $p$. The Schwartz kernel
$(xx')^{-\la+\ndemi}r(\la)$ of $\til{R}(\la)$ is meromorphic and
its polar part at $\la_0$ is the same as the one of
$(xx')^{-\la+\ndemi}(r_1(\la)+r_2(\la))$ since $r_0(\la)$ is
holomorphic in $\cc$. We then can easily check \cite[Prop.
3.3]{GU} that we have in $V_{\la_0}$
\beq\label{bjxi}
\Xi_{\la_0}(\til{R}(\la))=\sum_{j=-p}^{-1}B_j(\la-\la_0)^j
\eeq
where $B_j\in\mc{L}(x^{2\eps}L^2(X),x^{-2\eps}L^2(X))$ has a
Schwartz kernel of the form
\beq\label{noyaubj}
b_j(x,y,x',y')=\sum_{i=1}^{r_j} \psi_{ji}(x,y)\varphi_{ji}(x',y')
\left|\frac{dxdydx'dy'}{x^{n+1}x'^{n+1}}\right|^\demi, \quad
\psi_{ij},\varphi_{ij}\in x^{\ndemi}C^\infty(\bar{X}).
\eeq
Observe now that $x^{\la-\ndemi}$ has the following Taylor
expansion at $\la_0$
\[x^{\la-\ndemi}=x^{\la_0-\ndemi}\sum_{j=0}^{p-1}\log^j(x)\frac{(\la-\la_0)^j}{j!}+
O((\la-\la_0)^p)\]
in the sense of operators of $\mc{L}(x^NL^2(X),x^{2\eps}L^2(X))$ and
$\mc{L}(x^{-2\eps}L^2(X),x^{-N}L^2(X))$. We deduce that
$z'(\la)R(\la)$ has a residue $A$ satisfying
\[\textrm{Im}(A)\subset \sum_{j=0}^{p-1}x^{\la_0}\log^j(x) C^\infty(\bar{X},\Gamma_0^\demi)\]
and we are done.
\qed\\

\subsection{Scattering matrix}
Joshi and S\'a Barreto \cite{JSB} have shown that
the scattering matrix $S(\la)$ (defined in the introduction)
has the following Schwartz kernel
\beq\label{noyaudiffusion}
s(\la):=(2\la-n)\left (\beta_{\pl}\right )_*\left (\beta^*\left (
x^{-\la+\ndemi}x'^{-\la+\ndemi}r(\la)\right )|_{\mc{T}\cap\mc{B}}\right )
\eeq
Following (\ref{k2}) and (\ref{noyaudiffusion}) we have in
$\cc\setminus(\mc{R}\cup (\ndemi+\nn))$
\beq\label{noyaudesla}
s(\la)=\left (\beta_{\pl}\right )_*\left (r^{-2\la}k_1(\la)\right )+k_2(\la),
\eeq
\[k_1(\la)\in C^\infty(\pl\bar{X}\x_0\pl\bar{X},\Gamma^\demi), \quad
k_2(\la)\in C^\infty(\pl\bar{X}\x\pl\bar{X},\Gamma^\demi)\]
where $k_1(\la)$ and $k_2(\la)$ are meromorphic in $\la\in\cc$.
Outside its poles, $s(\la)$ is a conormal distribution of
order $-2\la$ associated to $\delta_{\pl\bar{X}}$
and $S(\la)$ is a pseudo-differential operator of order
$2\la-n$ on $\pl\bar{X}$.
In the sense of Shubin \cite[Def. 11.2]{SH}, $S(\la)$ is a
holomorphic family in
$\left\{\Re(\la)<\ndemi\right\}\setminus \mc{R}$
of zeroth order pseudo-differential operators.
We then deduce that $S(\la)$ is holomorphic in the same open set,
with values in $\mc{L}(L^2(\pl\bar{X}))$.
Recall the functional equation satisfied by $S(\la)$ (cf. \cite{GRZ})
\beq\label{eqfonct}
S(\la)^{-1}=S(n-\la)=S(\la)^* , \quad \Re(\la)=\ndemi,
\quad \la\not=\ndemi
\eeq
which also proves that $S(\la)$ is regular on the line
$\{\Re(\la)=\ndemi\}$.
Furthermore, (\ref{eqfonct}) holds also for $\til{S}(\la)$ and
by analytic extension we have on $\cc\setminus\mc{R}$
\[\til{S}^{-1}(\la)=\til{S}(n-\la).\]
The principal symbol of $S(\la)$ is given in (\ref{symboleprincipal})
and the renormalization $\til{S}(\la)$ of $S(\la)$ defined in
(\ref{stilde}) is Fredholm with index $0$, consequently we are
in the framework of Section 2.\\

Using Lemmas \ref{partiesinguliere} and (\ref{noyaudesla}),
we then obtain the
\begin{lem}\label{decompdesla}
Let $\la_0\in \{\Re(\la)<\ndemi\}$ a pole of $S(\la)$. Then
$\la_0\in\mc{R}$ and, following the notations of Lemma
\ref{partiesinguliere}, we have near $\la_0$
\beq\label{factorisationdesla}
S(\la)=(2\la-n)\trans\Phi^\sharp(\la)
F_1(\la)\left(\sum_{j=1}^m(z(\la)-z(\la_0))^{k_j}
P_j\right)F_2(\la)\Phi^\sharp(\la)+H^\sharp(\la)
\eeq
with
$H^\sharp(\la)\in\mc{H}ol(V_{\la_0},\mc{L}(L^2(\pl\bar{X})))$
and $\Phi^\sharp(\la)\in
\mc{H}ol(V_{\la_0},\mc{L}(L^2(\pl\bar{X}),\cc^q))$.
\end{lem}
\textsl{Proof}: the fact that $\la_0\in\mc{R}$ is
straightforward since if $r(\la)$ was holomorphic
one would have $s(\la)$ holomorphic in view of
(\ref{noyaudiffusion}). Now, $\til{R}(\la)$ being defined in
(\ref{rtilde}), we saw in Lemma \ref{partiesinguliere}
that the polar part of $\til{R}(\la)$ at $\la_0$
has a Schwartz kernel $\Xi_{\la_0}(\til{r}(\la))$ satisfying
\beq\label{structurexirtilde}
\Xi_{\la_0}(\til{r}(\la))\in
(xx')^\ndemi C^\infty(\bar{X}\x\bar{X},\Gamma_0^\demi).
\eeq
Let $\Phi(\la):=\sum_{i=0}^{p-1}\frac{(\la-\la_0)^i}{i!}
\frac{d^i}{d\la^i}(\Phi x^{-\la+\ndemi})|_{\la=\la_0}$ in the
sense of operators of $\mc{L}(x^{2\eps}L^2(X),\cc^q)$:
\[\Phi(\la) :\left\{
\begin{array}{rcl}
x^{2\eps}L^2(X) & \to & \cc^q\\
f & \to & \left(\sum_{j=0}^{p-1}\frac{(\la_0-\la)^j}{j!}\cjg
\log^j(x)x^{-\la_0+\ndemi}\psi_l,f\cjd\right)_{l=1,\dots,q}
\end{array}\right..\]
Lemma \ref{partiesinguliere} implies that \beq\label{contenue}
\Xi_{\la_0}(\til{R}(\la))=\Xi_{\la_0}\left(\trans\Phi(\la)F_1(\la)
\Big(\sum_{j=1}^m(z(\la)-z(\la_0))^{k_j}
P_j\Big)F_2(\la)\Phi(\la)\right). \eeq
Let $C:=\sum_{j=-p}^{-1}\textrm{Im}(B_j)$ with $B_j$ the operators
defined in (\ref{bjxi}) and let $\Pi_{C}$ be the orthogonal projection of
$x^{-2\eps}L^2(X)$ onto $C$. We multiply (\ref{contenue}) on the
left by $\Pi_C$ and on the right by $\trans \Pi_C$, and using
that $\Xi_{\la_0}(\til{R}(\la))$ is symmetric (since $\trans
R(\la)=R(\la)$) we deduce that (\ref{contenue}) remains true if
$\Phi(\la)$ is replaced by
\[\left\{
\begin{array}{rcl}
x^{2\eps}L^2(X) & \to & \cc^q\\
f & \to & \left(\sum_{j=0}^{p-1}\frac{(\la_0-\la)^j}{j!}
\cjg \Pi_C(\log^j(x)x^{-\la_0+\ndemi}\psi_l),
f\cjd\right)_{l=1,\dots,q}
\end{array}\right.\]
so that the logarithmic terms disappear.
Finally, we can use the representation of $S(\la)$ by its Schwartz
kernel (\ref{noyaudesla}) and we obtain
\[\Xi_{\la_0}(S(\la))=\Xi_{\la_0}\left((2\la-n)
\trans\Phi^\sharp(\la) F_1(\la)\Big
(\sum_{j=1}^m(z(\la)-z(\la_0))^{k_j}P_j\Big)F_2(\la)\Phi^\sharp(\la)\right ),\]
with
\[\Phi^\sharp(\la): \left\{\begin{array}{rcl}
L^2(\pl\bar{X}) & \to & \cc^q\\
f & \to & \left(\sum_{j=0}^{p-1}\frac{(\la_0-\la)^j}{j!}
\cjg \Pi_C(\log^j(x)x^{-\la_0+\ndemi}\psi_l)|_{\pl\bar{X}},
f\cjd\right)_{l=1,\dots,q}
\end{array}\right.,\]
the proof is achieved.
\qed\\

From this lemma, we deduce the
\begin{cor}\label{sens1}
If $\la_0\in \{\Re(\la)<\ndemi\}$ is
a pole of $S(\la)$, it is a pole of $R(\la)$ such that
\[m_{\la_0}(z'(\la)R(\la))\geq
N_{\la_0}\Big(c(n-\la)\til{S}(n-\la)\Big).\]
\end{cor}
\textsl{Proof}: firstly,
(\ref{factorisationdesla}) can be expressed by
\[c(\la)\til{S}(\la)=F_3(\la)\Big
(\sum_{j=1}^m(z(\la)-z(\la_0))^{k_j}
P_j\Big)F_4(\la)+\til{H}^\sharp(\la),\]
\[F_3(\la):=(2\la-n)\Lambda^{-\la+\ndemi}\trans\Phi^\sharp(\la) F_1(\la), \quad
F_4(\la):=F_2(\la)\Phi^\sharp(\la)\Lambda^{-\la+\ndemi},\]
\[\til{H}^\sharp(\la):=(2\la-n)\Lambda^{-\la+\ndemi}
H^\sharp(\la)\Lambda^{-\la+\ndemi}.\]
Note that we can take $k_1\leq\dots\leq k_m<0$
and set $(\varphi_0^{(j)})_{j=1,\dots,M}$
a canonical system of eigenvectors of
$c(n-\la)\til{S}(n-\la)$ at $\la_0$ with
$r_1\geq \dots\geq r_M$ the associated partial null
multiplicities (this canonical system exists
and is deduced from the one of $\til{S}(n-\la)$).
Let us show that $M\leq m$ and, by induction, that
$r_j\leq -k_j$ for all $j=1,\dots,M$.

If $\varphi^{(j)}(\la)$ is a root function
of $c(n-\la)\til{S}(n-\la)$ at $\la_0$ corresponding
to $\varphi_0^{(j)}$, there exists a holomorphic function $\phi^{(j)}(\la)$
such that
\[c(n-\la)\til{S}(n-\la)\varphi^{(j)}(\la)=
(z(\la)-z(\la_0))^{r_j}\phi^{(j)}(\la)\]
with $\phi^{(j)}(\la_0)\not=0$, hence when $\la$ approaches
$\la_0$ in the following identity
\[\varphi^{(j)}(\la)=\sum_{l=1}^m(z(\la)-z(\la_0))^{r_j+k_l}
F_3(\la)P_lF_4(\la)\phi^{(j)}(\la)+(z(\la)-z(\la_0))^{r_j}
\til{H}^\sharp(\la)\phi^{(j)}(\la),\]
we find that $r_1\leq-k_1$ and
$\varphi^{(j)}_0$ is in the vectorial space
\[E_j:=\textrm{Vect}\{F_3(\la_0)P_lF_4(\la_0)
L^2(\pl\bar{X});r_j\leq-k_l\}.\]
Moreover, the order on $(r_j)_{j=1,\dots,M}$ implies
that $E_j\subset E_M$ for $j=1,\dots,M$ but $\dim E_M\leq m$
since $\rang(P_l)=1$, thus we necessarily have $M\leq m$,
$(\varphi^{(j)}_0)_j$ being independent by assumption.
Now let $j\leq M$ and suppose that $r_i\leq -k_i$ for all
$i\leq j$. We first note that $E_j\subset E_{j+1}$ since
$r_{j+1}\leq r_j$. If $r_{j+1}>-k_{j+1}$, we would have
$\dim E_{j+1}\leq j$ but $E_{j+1}$ contains
the linearly independent vectors
$\varphi_0^{(1)},\dots,\varphi_0^{(j+1)}$, so a contradiction.
One concludes that $r_{j+1}\leq -k_{j+1}$ and
\[N_{\la_0}\Big(c(n-\la)\til{S}(n-\la)\Big)=
\sum_{j=1}^M{r_j}\leq -\sum_{l=1}^mk_l
=q=m_{\la_0}(z'(\la)R(\la)),\]
the corollary is proved.
\qed\\
\begin{lem}\label{sens2}
Let $\la_0\in \{\Re(\la)<\ndemi\}$ be a pole
of $R(\la)$ of finite multiplicity. If
$\la_0(n-\la_0)\notin\sigma_{pp}(\Delta_g)$ or $\la_0\notin \demi(n-\nn)$,
then $\la_0$ is a pole of $S(\la)$ such that
\[m_{\la_0}(z'(\la)R(\la))\leq
N_{\la_0}\Big(c(n-\la)\til{S}(n-\la)\Big).\]
\end{lem}
\textsl{Proof}: we first suppose that $\la_0$ is not a pole of
$c(\la)$ (i.e. $\la_0\notin \ndemi-\nn$). From Gohberg-Sigal theory,
one can factorize $\til{S}(\la)$ near $\la_0$ as in (\ref{factorisation})
\beq\label{fact1}
c(\la)\til{S}(\la)=U_1(\la)\left(
P_0+\sum_{l=1}^m(\la-\la_0)^{k_l}P_l\right)U_2(\la)
\eeq
with $U_1(\la)$, $U_2(\la)$ some holomorphically
invertible operators near $\la_0$ and
\[P_iP_j=\delta_{ij}P_j, \quad \rang (P_l)=1 \textrm{ for } l=1,\dots,m,
\quad 1=\sum_{j=0}^mP_j,\quad k_l\in\zz^*.\]
Take the Green equation between
the resolvent and the scattering operator
(see \cite{P,P1,G1,GZ3,GU})
\beq\label{mesurespect}
R(\la)-R(n-\la)=(2\la-n)\trans
E(n-\la)\Lambda^{\la-\ndemi}c(\la)\til{S}(\la)
\Lambda^{\la-\ndemi}E(n-\la)
\eeq
on $\mc{L}(x^NL^2(X),x^{-N}L^2(X))$ with $\ndemi-N<|\Re(\la)|<\ndemi$
and $E(\la)$ the transpose of the Eisenstein operator,
its Schwartz kernel being
\[e(\la):=\til{\beta}_*\left (\beta^*(x^{-\la+\ndemi}r(\la))
|_{\mc{T}}\right).\]
We can suppose that $k_1\leq\dots\leq k_m$ and set $p:=\max(0,-k_1)$. 
We consider the following Laurent expansions at $\la_0$
\beq\label{laurent1}
\begin{array}{c}
(n-2\la)R(n-\la)= \sum_{i=-1}^{p-1}
R_{i}(\la-\la_0)^i +O((\la-\la_0)^{p}),\\
(2\la-n)U_2(\la)\Lambda^{\la-\ndemi}E(n-\la)=\sum_{i=-1}^{p-1}
E^{(2)}_i(\la-\la_0)^i +O((\la-\la_0)^{p}),\\
(n-2\la)\trans E(n-\la)\Lambda^{\la-\ndemi}U_1(\la)= \sum_{i=-1}^{p-1}
E^{(1)}_i(\la-\la_0)^i +O((\la-\la_0)^{p}),
\end{array}
\eeq
where $R_{-1}$ and $E^{(j)}_{-1}$ are not $0$ if and only if
$\la_0(n-\la_0)\in\sigma_{pp}(\Delta_g)$, and in this case
\beq\label{structre}\begin{array}{c}
R_{-1}=-\sum_{i=1}^k
\phi_i\otimes\phi_i,\\
E^{(2)}_{-1}=\sum_{i=1}^k
U_2(\la_0)\Lambda^{\la_0-\ndemi}(x^{\la_0-\ndemi}\phi_i)
|_{\pl\bar{X}}\otimes\phi_i,\\
E^{(1)}_{-1}= -\sum_{i=1}^k\phi_i\otimes \trans
U_1(\la_0)\Lambda^{\la_0-\ndemi}(x^{\la_0-\ndemi}\phi_i)
|_{\pl\bar{X}},\end{array} \eeq with $\phi_i\in
x^{n-\la_0}C^\infty(\bar{X},\Gamma_0^\demi)$ the normalized
eigenfunctions of $\Delta_g$ for the eigenvalue $\la_0(n-\la_0)$.
From (\ref{fact1}), (\ref{mesurespect}) and (\ref{laurent1}) we
obtain
\beq\label{resrla1}
A:=\textrm{Res}_{\la_0}((n-2\la)R(\la))
=R_{-1} +\sum_{\substack{j+i+k_l=-1\\k_l\geq 0}}
E^{(1)}_iP_lE^{(2)}_j +\sum_{\substack{j+i+k_l=-1\\k_l<0}}
E^{(1)}_iP_lE^{(2)}_j \eeq
where by convention $k_l=0\iff l=0$. We set $V:=\textrm{Im}(A_1)+\textrm{Im}(A_2)$ with
\[A_1:=R_{-1}+E^{(1)}_{-1}P_0E^{(2)}_0+
E^{(1)}_{-1}\left(\sum_{k_l=1}
P_l\right)E^{(2)}_{-1},\]
\[A_2:=E^{(1)}_{0}P_0E^{(2)}_{-1}+\sum_{\substack{j+i+k_l=-1\\k_l<0}} E^{(1)}_iP_lE^{(2)}_j.\]
Remark from (\ref{structre}) that
\[\textrm{Im}(A_1)\subset x^{n-\la_0}C^\infty(\bar{X},\Gamma_0^\demi),
\quad (\Delta_g-\la_0(n-\la_0))A_1=0\]
and in view of Lemma \ref{partiesinguliere} we know that there exists
$p\in\nn$ such that
\[ \textrm{Im}(A)\subset
\sum_{j=0}^{p-1}x^{\la_0}\log^j(x)C^\infty(\bar{X},\Gamma_0^\demi),\quad
(\Delta_g-\la_0(n-\la_0))^pA=0\]
thus we can argue that
\[\forall u\in V,\quad (\Delta_g-\la_0(n-\la_0))^pu=0. \]
Note that if $\la_0\notin \demi(n-\nn)$, we clearly have
\[x^{n-\la_0}C^\infty(\bar{X},\Gamma_0^\demi)\cap
\sum_{j=0}^{p-1}x^{\la_0}\log^j(x)C^\infty(\bar{X},\Gamma_0^\demi)
\subset \dot{C}^\infty(\bar{X},\Gamma^\demi_0),\]
therefore, if $V_1,V_2$ are defined by
\[V_1=V\cap x^{n-\la_0}C^\infty(\bar{X},\Gamma^\demi_0),\quad
V_2=V\cap \sum_{j=0}^{p-1}x^{\la_0}\log^j(x) C^\infty(\bar{X},\Gamma_0^\demi),\]
we deduce from the unique continuation principle proved by Mazzeo \cite{MA} that
\[V_1\cap V_2\subset \dot{C}^\infty(\bar{X},\Gamma^\demi_0)\cap
\ker_{L^2}(\Delta_g-\la_0(n-\la_0))^p=0.\]
Hence, we can split $V=V_1\oplus V_2\oplus V_3$ with $V_3$ a direct complement
of $V_1\oplus V_2$ in $V$. Let $\Pi_{V_2}$ be the projection of $V$ onto
$V_2$ parallel to $V_1\oplus V_3$, $\Pi_V$ the orthogonal projection of $x^{-N}L^2(X)$ onto
$V$ and $\iota_V$ the inclusion of $V$ into $x^{-N}L^2(X)$. We
multiply (\ref{resrla1}) on the left by $\Pi'_{V_2}:=\iota_V\Pi_{V_2}\Pi_V$
and on the right by $\trans \Pi'_{V_2}$ to obtain
\[A=\sum_{\substack{j+i+k_l=-1\\k_l<0}}\Pi'_{V_2}
E^{(1)}_iP_lE^{(2)}_j\trans\Pi'_{V_2}\]
by construction of $V_2$ and using the symmetry $\trans A=A$ (since
$\trans R(\la)=R(\la)$).
Now remark that
\[\sum_{\substack{j+i+k_l=-1\\k_l<0}}\Pi'_{V_2}
E^{(1)}_iP_lE^{(2)}_j\trans\Pi'_{V_2}=
\sum_{k_l<0}\sum_{i=0}^{-k_l-1}\Pi'_{V_2}
 E^{(1)}_iP_lE^{(2)}_{-k_l-1-i}\trans\Pi'_{V_2}\]
and the rank of this operator is bounded by
$-\sum_{k_l<0}k_l=N_{\la_0}(c(n-\la)\til{S}(n-\la))$
since $\rang(P_l)=1$.
The lemma is proved when $\la_0\notin \ndemi-\nn$.\\

On the other hand if $\la_0\in \ndemi-\nn$ and $\la_0(n-\la_0)\notin
\sigma_{pp}(\Delta_g)$, we have $R_{-1}=0$, $E^{(1)}_{-1}=0$
and $E^{(2)}_{-1}=0$ in (\ref{laurent1}). Therefore,
the same proof works if we replace
(\ref{fact1}) and (\ref{resrla1}) by
\[c(\la)\til{S}(\la)=U_1(\la)\left((\la-\la_0)P_0+
\sum_{l=1}^m(\la-\la_0)^{k_l+1}P_l\right)U_2(\la),\]
\[\textrm{Res}_{\la_0}((n-2\la)R(\la))
=\sum_{\substack{j+i+k_l=-2\\ k_l<-1}}E^{(1)}_iP_lE^{(2)}_j\]
the first one being obtained from Gohberg-Sigal factorization
(\ref{factorisation}) of $\til{S}(\la)$ at $\la_0$.
Now observe that the rank of
\[\sum_{\substack{j+i+k_l=-2 \\ k_l<-1}}
\Pi'_{V_2}E^{(1)}_iP_lE^{(2)}_j\trans\Pi'_{V_2}=
\sum_{k_l<-1}\sum_{i=0}^{-k_l-2}
\Pi'_{V_2}E^{(1)}_iP_lE^{(2)}_{-k_l-2-i}\trans\Pi'_{V_2}\]
is bounded by
\[-\sum_{k_l<-1}(k_l+1)=-\sum_{k_l<0}(k_l+1)=
N_{\la_0}(\til{S}(n-\la))-\dim\ker_{\la_0} \til{S}(n-\la)=
N_{\la_0}\Big(c(n-\la)\til{S}(n-\la)\Big)\]
in view of (\ref{ndela0}), the proof is complete.
\qed\\

\textsl{Proof of Theorem \ref{multiplicite}}:
we combine the results of Corollary \ref{sens1} and Lemma
\ref{sens2} with (\ref{ndela0}) and (\ref{gohbergsigal}),
and observe that 
\[\ker_{\la_0}\til{S}(n-\la)=\ker
\til{S}(n-\la_0)=\ker\textrm{Res}_{n-\la_0}S(\la),\]
then it remains to show that
\beq\label{ndestilde}
N_{\la_0}(\til{S}(\la))=m_{n-\la_0}.
\eeq
Whereas the case
$\la_0(n-\la_0)\notin \sigma_{pp}(\Delta_g)$
is clear since $\til{S}(\la)^{-1}=\til{S}(n-\la)$ is holomorphic
near $\la_0$ and $m_{n-\la_0}=0$, the case
$\la_0(n-\la_0)\in\sigma_{pp}(\Delta_g)$ needs a little
more care.
In view of (\ref{residuvp}) and (\ref{noyaudiffusion}),
$\til{S}(\la)$ has the following polar part at $n-\la_0$
\[C(\la_0)(\la-n+\la_0)^{-1}\sum_{j=1}^k \Lambda^{\la_0-\ndemi}\phi_j^\sharp\otimes
\Lambda^{\la_0-\ndemi}\phi_j^\sharp\]
with $C(\la_0)\not=0$ if $\la_0\notin \ndemi-\nn$, $k=m_{n-\la_0}$ and
$\phi_j^\sharp:=x^{\la_0-\ndemi}\phi_j|_{\pl\bar{X}}$ (where $(\phi_j)_j$
is an orthonormal basis of $\ker_{L^2}(\Delta_g-\la_0(n-\la_0))$ 
as in (\ref{residuvp})).
It is not difficult to see that $(\phi_j^\sharp)_j$ are independent,
otherwise there would exist a non zero solution 
$u\in x^{n-\la_0+1}C^\infty(\bar{X},\Gamma^\demi_0)$
of $(\Delta_g-\la_0(n-\la_0))u=0$ and a Taylor expansion
of this equation at $x=0$ proves that
$u\in\dot{C}^\infty(\Bar{X},\Gamma^\demi_0)$,
which is excluded according to Mazzeo's result \cite{MA}. 
Since the pole is a first order pole, the factorization
of $\til{S}(\la)$ as in (\ref{factorisation})
near $n-\la_0$ is clear for the $k_l<0$:
we have $m=k$ and $k_l=-1$ for $l=1,\dots,k$. 
Using (\ref{factorisation2}) and $\til{S}(\la)^{-1}=\til{S}(n-\la)$, 
one then obtain that the partial null multiplicities 
of $\til{S}(\la)$ at $\la_0$ are $\{-k_1,\dots,-k_k\}$ 
which gives (\ref{ndestilde}) and the theorem.
\qed\\

\end{document}